\documentstyle[11pt]{article}
\newtheorem{thm}{Theorem}
\newtheorem{rk}{Remark}
\newtheorem{prop}{Proposition}

\begin{document}

\centerline{\large{\bf\underline{"Vector bundles" over quantum Heisenberg
manifolds.}}}
\vspace{.2in}
{\centerline {\large{Beatriz Abadie\footnote {The material of this work is part
of the author's Ph.D dissertation, submitted to the University of California at
Berkeley in May 1992.}}}}
\vspace{.2in}
{\small \underline{{\bf Abstract}}:  We construct, out of Rieffel projections,
projections in certain algebras which are strong-Morita equivalent to the
quantum Heisenberg manifolds $D^{c}_{\mu\nu}$. Then, by means of techniques
from the Morita equivalence theory, we get finitely generated and projective
modules over the algebras $D_{\mu\nu}^{c}$. This enables us to show that the
group $Z+2\mu Z +2\nu Z$ is contained in the range of the trace on
$K_{0}(D^{c}_{\mu\nu})$.}
\vspace{.3in}
\newline\underline{{\bf Preliminaries.}}
Let $G$ be the Heisenberg group,
\[G=\left\{\left( \begin{array}{ccc}
        1 & y & z \\
        0 &1 & x  \\
        0 &0 &1
        \end{array}
        \right) : x,y,z\in R \right\}, \]
and, for a positive integer $c$, let $H_{c}$ be the subgroup of $G$ obtained
when $x$, $y$, and $cz$ are integers. The Heisenberg manifold $M_{c}$ is the
quotient $G/H_{c}$.

Non-zero Poisson brackets on $M_{c}$ that are invariant under the action $G$ on
$M_{c}$ by left translation can be parametrized by two real numbers $\mu$ and
$\nu$, with $\mu^{2}+\nu^{2}\not = 0$ (\cite{rfhm}).

For each positive integer $c$ and real numbers $\mu$ and $\nu$ as above,
Rieffel constructed (\cite{rfhm}) a deformation quantization
$\{D^{c,\hbar}_{\mu\nu}\}_{\hbar\in R}$ of $M_{c}$ in the direction of the
Poisson bracket $\Lambda_{\mu\nu}$.

Since $D_{\mu\nu}^{c,\hbar}$ is isomorphic to $D^{c,1}_{\hbar\mu,\hbar\nu}$,
and we will not need to keep track of the Planck constant $\hbar$, we absorb it
from now on into the parameters $\mu$ and $\nu$. Thus we will use
$D_{\mu\nu}^{c}$ to denote $D^{c,1}_{\mu\nu}$.

As shown in \cite{rfhm}, the algebra $D_{\mu\nu}^{c}$ can be described as the
generalized fixed-point algebra of the crossed-product $C_{0}(R\times
T)\times_{\lambda}Z$, where $\lambda_{p}(x,y)=(x+2p\mu,y+2p\nu)$, for all $p\in
Z$,  under the action $\rho$ of $Z$ defined by
\[ (\rho_{k}\Phi)(x,y,p)=e(ckp(y-p\nu))\Phi(x+k,y,p),\]
where $k$, $p\in Z$, $\Phi\in C_{c}(R\times T\times Z)$, and, for any real
number $x$, $e(x)=exp(2\pi i x)$.

The action $\rho$ defined above corresponds to the action $\rho$ defined in
\cite[p.539]{rfhm}, after taking Fourier transform in the third variable to get
the algebra denoted in that paper by $A_{\hbar}$, and viewing $A_{\hbar}$ as a
dense *-subalgebra of $C_{0}(R\times T)\times_{\lambda}Z$ via the embedding $J$
defined in \cite[p.547]{rfhm}.
Equivalently, $D^{c}_{\mu\nu}$ is the closure in the multiplier algebra of
$C_{0}(R\times T)\times_{\lambda}Z$ of the *-subalgebra $D_{0}$ consisting of
functions $\Phi\in C(R\times T\times Z)$ which have compact support on $Z$ and
satisfy
\[ \Phi(x+k,y,p)=e(-ckp(y-p\nu))\Phi(x,y,p),\]
for all $k,p\in Z$, and $(x,y)\in R\times T$ ($D_{0}$ is the image under the
embedding $J$ mentioned above of the subalgebra denoted by $C^{\rho}$ in the
proof of \cite[Thm.5.4]{rfhm}).

There is a faithful trace (\cite{rfhm}) $\tau_{D}$ on $D_{\mu\nu}^{c}$ defined
for $\Phi\in D_{0}$, by
\[ \tau_{D}(\Phi)=\int_{T^{2}}\Phi(x,y,0)dxdy.\]
It can be shown (\cite{fpa}) that the algebra $D^{c}_{\mu\nu}$ is strong-Morita
equivalent to the generalized fixed-point algebra $E_{\mu\nu}^{c}$ of the
crossed product $C_{0}(R\times T)\times_{\sigma} Z$ under the action $\gamma$
of $Z$, where $\sigma_{k}(x,y)=(x-k,y)$ and
\[(\gamma_{p}\Phi)(x,y,k)=e(-ckp(y-p\nu))\Phi(x-2p\mu,y-2p\nu,k),\]
for  $k$, $p\in Z$ and  $\Phi\in C_{c}(R\times T\times Z)$.

As for the quantum Heisenberg manifolds case, $E_{\mu\nu}^{c}$ can also be
described (see \cite{fpa}) as the closure in the multiplier algebra of
$C_{0}(R\times T)\times_{\sigma}Z$ of the *-algebra $E_{0}$ consisting of
functions $\Phi\in C(R\times T\times Z)$, with compact support on Z and such
that
\[\Phi(x-2p\mu,y-2p\nu,k)=e(ckp(y-p\nu))\Phi(x,y,k),\]
for all $k,p \in Z$, $(x,y)\in R\times T$.
The equivalence $D^{c}_{\mu\nu}$-$E_{\mu\nu}^{c}$ bimodule $X$ constructed in
\cite{fpa} is the completion of $C_{c}(R\times T)$ with respect to either one
of the norms induced by the $D_{\mu\nu}^{c}$ and  $E_{\mu\nu}^{c}$-valued inner
products given by
\[<f,g>_{D}(x,y,p)=\sum_{k\in
Z}e(ckp(y-p\nu))f(x+k,y)\overline{g}(x-2p\mu+k,y-2p\nu)\]
and
\[<f,g>_{E}(x,y,k)=\sum_{p\in
Z}e(-ckp(y-p\nu))\overline{f}(x-2p\mu,y-2p\nu)g(x-2p\mu+k,y-2p\nu),\]
respectively, where $f,g\in C_{0}(R\times T)$, $\Phi\in D_{0}$, $\Psi\in
E_{0}$, $(x,y)\in R\times T$, and $k,p\in Z$.

In what follows we produce finitely generated and projective modules over the
algebras $D_{\mu\nu}^{c}$. To do this we apply to the Morita equivalence
structure described above the methods for constructing projections provided by
the Morita equivalence theory. Finally, we get a partial description of the
range of the trace at the level of $K_{0}(D_{\mu\nu}^{c})$.

\begin{rk}
\label{morita}
 First notice that both $D_{0}$ and $E_{0}$ have identity elements $I_{D}$ and
$I_{E}$, respectively, defined by
\[I_{D}(x,y,p)=\delta_{0}(p)\mbox{ and } I_{E}(x,y,k)=\delta_{0}(k),\]
for $(x,y)\in R\times T$ and $k,p\in Z$.

Therefore, by \cite[Prop. 1.2]{rfhd}, if $P$ is a projection in $E_{0}$, then
$XP$ is a projective finitely
generated left module over $D_{\mu\nu}^{c}$, and the corresponding
projection in $M_{m}(D_{\mu\nu}^{c})$ is given by

 \[Q =\left( \begin{array}{ccc}
           <y_{1},x_{1}>_{D} & ... &<y_{m},x_{1}>_{D}\\
           ...   & ... & ... \\
           <y_{1},x_{m}>_{D} & ... & <y_{m},x_{m}>_{D}
           \end{array}
           \right) \]
\vspace{.1in}
where, for $i=1,...,m$ ,  $x_{i},y_{i}\in XP$ are such that
$P=\sum_{i=1}^{i=m}<x_{i},y_{i}>_{E}$.

On the other hand (\cite[Prop. 2.2]{rfirr}), the trace $\tau_{D}$ on
$D^{c}_{\mu\nu}$ induces a trace $\tau_{E}$ on $E_{\mu\nu}^{c}$ via
\[\tau_{E}(<f,g>_{E})=\tau_{D}(<g,f>_{D}).\]
A straightforward computation shows that for $\Psi\in E_{0}$ we have
\[ \tau_{E}(\Psi)=\int_{0}^{2\mu}\int_{0}^{1}\Psi(x,y,0)dxdy.\]
Then, in the notation above we get
\[\tau_{D}(Q)=\sum_{i=1}^{1=m}\tau_{D}(<y_{i},x_{i}>_{D})=\sum_{i=1}^{i=m}\tau_{E}(<x_{i},y_{i}>_{E})=\tau_{E}(P).\]
\end{rk}
\vspace{.1in}

\begin{thm}
\label{mu}
The bimodule $X$ is a finitely generated and projective $D^{c}_{\mu\nu}$-module
of trace $2\mu$. If $\nu\in [0,1/2]$, and $\mu>1$, then there is a finitely
generated projective $D_{\mu\nu}^{c}$-submodule of $X$ with trace $2\nu$.
\end{thm}
\vspace{.1in}
\underline{\em{Proof}}:

\vspace{.1in}
Let us take $P=I_{E}$, in the notation of Remark \ref{morita}. Then $X=XP$ is
finitely generated and projective and its trace is $\tau_{E}(I_{E})=2\mu$.

We now find a projection $P$ in $E_{0}$ with $\tau_{E}(P)=2\nu$, when $\nu\in
[0,1/2]$ and $\mu>1$, which ends the proof, in view of Remark \ref{morita}.
\newline So let us consider self-adjoint elements P of the form:
\[P(x,y,p)=f(x,y)\delta_{1}(p)+h(x,y)\delta_{0}(p)+\overline{f}(x-1,y)\delta_{-1}(p),\]
where $h$ and $f$ are bounded functions on $R\times T$ and $h$ is real-valued.
Our next step is to get functions $f$ and $h$ such that $P$ is a projection in
$E_{\mu\nu}^{c}$.

Now,
\[(P*P)(x,y,p)=\sum_{q\in Z}P(x,y,q)P(x+q,y,p-q),\]
and it follows that $P*P=P$ if and only if, for all $(x,y)\in R\times T$:
\vspace{.1in}

1) $f(x,y)f(x+1,y)=0$\\

2) $f(x,y)[h(x+1,y)+h(x,y)-1]=0$\\

3) $|f(x,y)|^{2}+|f(x-1,y)|^{2}=h(x,y)(1-h(x,y)).$
\vspace{.1in}
\newline We also want $P$ to be in $E_{0}$, so we require
\[P(x,y,p)=e(cp(y+\nu))P(x+2\mu,y+2\nu,p)\mbox{, that is}\]

4) $f(x,y)=e(c(y+\nu))f(x+2\mu,y+2\nu)$\\
and

5) $h(x,y)=h(x+2\mu,y+2\nu).$\\

It was shown on \cite[1.1]{rfirr} that for any $\zeta\in [0,1/2]$ there are
maps $F,H\in C(T)$ such
that:

1)' $F(t)F(t-\zeta)=0$\\

2)' $F(t)[1-H(t)-H(t-\zeta)]=0$\\

3)' $H(t)[1-H(t)]=|F(t)|^{2}+|F(t+\zeta)|^{2}$\\

4)' $\int_{T}H=\zeta$\\

5)' $0\leq H(t)\leq 1$ for any $t\in T$ and F vanishes on $[1/2,1]$.\\

Let us assume that $\nu\in [0,1/2]$, $\mu>1$ and let F and H be
functions satisfying 1)'-5)' for $\zeta=\nu/\mu$.

Translation of $t$ by $\zeta$ in equations 1)'-5)' plays the same role as
translation of $x$ by $1$ in equations 1)-5), which suggests taking $\zeta x$
as the variable $t$. However, the variable $y$ will play an important role in
getting $f$ and $h$ to  satisfy 4) and 5), for which we need to take
$t=1/2+y-\zeta x$.
\newline So let
\vspace{.05in}
\[h(x,y)=H(1/2+y-\zeta x),\] so h is in $C(R\times T)$, and it is real-valued
and bounded.
\newline Also,
\[h(x+2\mu,y+2\nu)=H(1/2+y+2\nu-\zeta x -2\nu)=H(1/2+y-\zeta x)=h(x,y),\]
so h satisfies 5).
\newline Now, for $(x,y)\in [0,2\mu]\times[0,1]$,  set

\[f(x,y)= \left\{ \begin{array}{ll}
                   F(1/2+y-\zeta x) & \mbox{if $y\leq x/(2\mu) $}\\
                  e(c(y+\nu))F(1/2+y-\zeta x ) & \mbox{if $y\geq
x/(2\mu)$}
                   \end{array}
                   \right. \]

To show f is continuous it suffices to prove that $F(1/2+y-\zeta x)=0$ when
$y=x/(2\mu)$, and that follows from the fact that $F$ vanishes on $[1/2,1]$,
and from the conditions on $\mu$ and $\nu$.
\newline Now, since $f(x,1)=f(x,0)$, f is continuous  on $[0,2\mu]\times
T$. We want to extend $f$\ to $R\times T$ by letting
\[f(x+2\mu,y)=e(-c(y-\nu))f(x,y-2\nu), \]
so as to have $f$ satisfy 4).
We only need to show that
\[f(2\mu,y)=e(-c(y-\nu))f(0,y-2\nu) \mbox{ for any } y\in T.\]
For an arbitrary $y\in R$, let $k,l\in Z$ be such that
$y+k,y-2\nu+l\in [0,1]$.
Then,
\[f(2\mu,y)=F(1/2+y+k-2\nu)=F(1/2+y-2\nu),\mbox{  and }\]
\[f(0,y-2\nu)=f(0,y-2\nu+l)=e(c(y-\nu+l))F(1/2+y-2\nu)=\]
\[=e(c(y-\nu))f(2\mu,y),\]
as wanted, and f, extended to $R\times T$ as above, satisfies 4). It remains to
show that f and g satisfy 1), 2) and 3):

\vspace{.1in}
1) $|f(x,y)f(x+1,y)|=|F(1/2+y-\zeta x)F(1/2+y-\zeta x-\zeta|=0$, by
1)'.\\

2) $|f(x,y)[h(x+1,y)+h(x,y)-1]|=|F(1/2+y-\zeta x)[H(1/2+y-\zeta
x-\zeta)+H(1/2+y-\zeta x)-1]|=0$, by 2)'.\\

3) $|f(x,y)|^{2}+|f(x-1,y)|^{2}=|F(1/2+y-\zeta x )|^{2}+|F(1/2+y-\zeta x
+\zeta)|^{2}=$
\newline $=H(1/2+y-\zeta x)[1-H(1/2+y-\zeta
x)]=h(x,y)(1-h(x,y))$,\hspace{.05in} by 3)'.
\vspace{.1in}

Therefore P is a projection on $E_{0}$. Besides,
\[\tau_{E}(P)=\int_{0}^{2\mu}\int_{T}h(x,y)dydx=\]
\[\int_{0}^{2\mu}(\int_{T}H(1/2 +y -\zeta
x)dy)dx=\int_{0}^{2\mu}\zeta=2\mu\zeta=2\nu,\mbox{ by 5)'.}\]

\begin{flushright}
Q.E.D.
\end{flushright}
\vspace{.1in}
The following propositions enable us to extend the previous results by dropping
the restrictions on $\mu$ and $\nu$.
\vspace{.1in}
\newline {\bf Notation:} In Propositions \ref{plusk} and \ref{minus} $\Pi$
denotes the faithful representation of $D_{\mu\nu}^{c}$ on $L^{2
}(R\times T\times Z)$ obtained by restriction of the left regular
representation of the multiplier algebra of $C_{0}(R\times T)\times_{\lambda}Z$
on $L^{2}(R\times T\times Z)$, i.e.
\[(\Pi_{\Phi}\xi)(x,y,p)=\sum_{q\in Z}\Phi(x+2p\mu,y+2p\nu,q)\xi(x,y,p-q),\]
for $\Phi\in D_{0}$, $\xi\in L^{2}(R\times T\times Z)$, and $(x,y,p)\in R\times
T\times Z$.

Notice that $\Pi$ is faithful because $Z$ is amenable (\cite[7.7.5 and
7.7.7.]{pd}).

\begin{prop}
\label{plusk} There is a trace-preserving isomorphism between $D_{\mu\nu}^{c}$
and $D_{\mu+k,\nu+l}^{c}$, for all $k,l\in Z$.
\end{prop}

\underline{\em{Proof}}:
\vspace{.05in}

It is clear that $\Phi\mapsto\Phi$ is an isomorphism between $D_{\mu\nu}^{c} $
and $D_{\mu,\nu+l}^{c}$, so let us assume $l=0, k=1.$
\newline Let $J:D_{\mu+1,\nu}^{c}\longrightarrow D_{\mu\nu}^{c}$ be defined at
the level of functions in $D_{0}$ by:
\[(J\Phi)(x,y,p)=e(c(4p^{3}\nu/3-p^{2}y))\Phi(x,y,p).\]
It is easily checked that $J\Phi\in D_{\mu\nu}^{c}$ for all $\Phi\in
D_{\mu+1,\nu}^{c}$.
Besides, the unitary operator  $U:L^{2}(R\times T\times Z)\longrightarrow
L^{2}(R\times T\times Z)$ given by
\[U\xi(x,y,p)=e(c(-4p^{3}\nu/3-p^{2}y))\xi(x,y,p)\]
intertwines $\Pi_{J\Phi}$ and $\Pi_{\Phi}$:
\[(\Pi_{\Phi}U\xi)(x,y,p)=\sum_{q\in
Z}\Phi(x+2p(\mu+1),y+2p\nu,q)U\xi(x,y,p-q)=\]
\[=\sum_{q\in Z}e(-2pcq(y+(2p-q)\nu)e(c((-4\nu/3)(p-q)^{3}-(p-q)^{2}y))\cdot\]
\[\cdot\Phi(x+2p\mu,y+2p\nu,q)\xi(x,y,p-q)=\]
\[=e(c(-4\nu p^{3}/3-p^{2}y))\sum_{q\in
Z}e(c(4q^{3}\nu/3-q^{2}(y+2p\nu))\Phi(x+2p\mu,y+2p\nu,q)\xi(x,y,p-q)=\]
\[=(U\Pi_{J\Phi}\xi)(x,y,p).\]

Also,
\[\tau(J\Phi)=\int_{0}^{1}\int_{T}J\Phi(x,y,0)=\int_{0}^{1}\Phi(x,y,0)=\tau(\Phi).\]
\begin{flushright}
Q.E.D.
\end{flushright}
\begin{prop}
\label{minus} There is a trace-preserving isomorphism between $D_{\mu\nu}^{c}$
and $D_{-\mu ,-\nu}^{c}$.
\end{prop}
\vspace{.1in}
\underline{\em{Proof}}:
\vspace{.05in}

Let $J:D_{\mu\nu}^{c}\longrightarrow D_{-\mu,-\nu}^{c}$ be defined, at the
level of functions, by:
\[(J\Phi)(x,y,p)=\Phi(-x,-y,p).\]

It is easily checked that $J\Phi\in D_{-\mu,-\nu}$. Besides, the unitary
operator
\newline $U:L^{2}(R\times T\times Z)\longrightarrow L^{2}(R\times T\times Z)$
defined by
\[(U\xi)(x,y,p)=\xi(-x,-y,p)\]
intertwines $\Pi_{\Phi}$ and $\Pi_{J\Phi}$:
\[[\Pi_{J\Phi}(U\xi)](x,y,p)=\sum_{q\in
Z}(J\Phi)(x-2p\mu,y-2p\nu,q)\xi(-x,-y,p-q)=\]
\[ =\sum_{q\in Z}\Phi(-x+2p\mu,-y+2p\nu,q)\xi(-x,-y,p-q)=\]
\[=(\Pi_{\Phi}\xi)(-x,-y,p)=(U\Pi_{\Phi}\xi)(x,y,p).\]

Finally, $J$ preserves the trace:
\vspace{.05in}
\[\tau(J\Phi)=\int_{T^{2}}\Phi(-x,-y,0)
=\tau(\Phi).\]
\begin{flushright}
Q.E.D.
\end{flushright}
\vspace{.1in}
\begin{thm} The range of the trace on $K_{0}(D_{\mu\nu}^{c})$ contains the set
$Z+2\mu Z+2\nu Z$.
\end{thm}
\underline{\em{Proof}}:

\vspace{.05in}
We obviously have $Z\subseteq \tau(K_{0}(D_{\mu\nu}^{c}))$, since
$D_{\mu\nu}^{c}$ has an identity element.
Besides, it follows from Theorem \ref{mu} that $2\mu Z\subseteq
\tau(K_{0}(D_{\mu\nu}^{c}))$.
So it only remains to show that $2\nu Z\subseteq \tau(K_{0}(D_{\mu\nu}^{c}))$.

Let $k\in Z$ be such that $\nu'=\pm\nu+k$ and $\nu'\in [0,1/2]$. Then one can
find $l\in Z$ and $\mu'=\pm \mu+l$ such that $\mu'\geq 1$. Thus , by
Propositions \ref{plusk} and \ref{minus} we have that
$\tau(K_{0}(D_{\mu'\nu'}^{c}))=\tau(K_{0}(D_{\mu\nu}^{c}))$.

Now, by Theorem \ref{mu} there is a projection in $M_{m}(D_{\mu'\nu'}^{c})$ for
some positive integer $m$ with trace $2\nu'=\pm 2\nu+2k$, which ends the proof.

\begin{flushright}
Q.E.D.
\end{flushright}

{\bf Remark.} It can  be shown (\cite{th}) that the inclusion in the previous
theorem is actually an identity.
\vspace{.1in}

{\bf Acnowledgement.} I am glad to thank my thesis advisor, Marc Rieffel for
his constant support and for many helpful suggestions and comments.

\vspace{.3in}
\begin{flushright}
Centro de Matem\'atica\\
Eduardo Acevedo 1139\\
CP 11200\\
Montevideo-Uruguay.\\
e-mail: abadie@cmat.edu.uy\\
\end{flushright}
\end{document}